\documentclass{amsart}
\usepackage{amsmath}
\usepackage{amssymb}
\usepackage{amscd}
\usepackage{amsthm}
\input{tableau}

\newcommand{\C}{{\mathbb C}}

\newcommand{\ngo}{\Omega}
\newcommand{\wts}{\Lambda^+}

\newcommand{\tD}{D}
\newcommand{\sD}{\tD^\circ}
\newcommand{\tq}{{\tilde q}}
\newcommand{\prho}{E}
\newcommand{\prj}{\kappa}

\newcommand{\tpi}{\tau}

\newcommand{\cE}{\mathcal{E}}

\newcommand{\pp}[1]{{\mathbf{p_{#1}}}}
\newcommand{\PP}[1]{{\mathbf{P_{#1}}}}
\newcommand{\Ar}{\text{${\mathbf A}$}}
\newcommand{\Al}{\text{${\mathbf A}^{-1}$}}
\newcommand{\Br}{\text{${\mathbf B}$}}
\newcommand{\Bl}{\text{${\mathbf B}^{-1}$}}
\newcommand{\Cr}{\text{${\mathbf C}$}}

\renewcommand{\tabstyle}{\scriptstyle}

\edef\det{\det\nolimits}

\theoremstyle{plain}

\newtheorem*{prop78}{Proposition~7.8}
\newtheorem*{prop79}{Proposition~7.9}
\newtheorem*{prop83}{Proposition~8.3}
\newtheorem*{prop85}{Proposition~8.5}
\newtheorem*{cor84}{Corllary~8.4}

\theoremstyle{definition}

\theoremstyle{remark}

\newcommand{\inv}{^{-1}}

\makeatletter
\newcounter{casify}
  {\begin{list}{{\it Case \arabic{casify}}. }
    {\usecounter{casify}
      \setlength{\leftmargin}{0pt}
      \setlength{\rightmargin}{0pt}
      \setlength{\labelwidth}{0pt}
      \setlength{\labelsep}{0pt}
      \setlength{\itemindent}{\parindent}}}
  {\end{list}}

\let\enumdepth\@enumdepth
\newcounter{enumsave}
\newenvironment{enumgen}[1]
  {\setcounter{enumsave}{\enumdepth}
     \enumdepth=#1
     \begin{enumerate}
     \ifnum\enumdepth>1  \fi
     \ifnum\enumdepth>2  \fi
     \ifnum\enumdepth>3  \fi
     \ifnum\enumdepth>4  \fi}
  {\end{enumerate} \enumdepth=\theenumsave}
\newcounter{enumpausectr}
\def\enumpause#1{\setcounter{enumpausectr}{\value{\@enumctr}}\end{enumerate}
  #1\begin{enumerate}\setcounter{\@enumctr}{\value{enumpausectr}}}
\makeatother

\title[Equivariant \emph{K}-theory in the general linear group]{Corrigendum to ``On the
equivariant \emph{K}-theory of the nilpotent cone in the general
linear group''}

\author{Pramod N. Achar}
\address{Department of Mathematics\\
  Louisiana State University\\
  Baton Rouge, LA 70803}
\email{pramod@math.lsu.edu}

\thanks{The author was partially supported NSF Grant No.~DMS-1500890.}

\begin{document}

\begin{abstract}
In the paper [P.~Achar, \emph{On the equivariant $K$-theory of the nilpotent cone in the general linear group}, Represent.~Theory {\bf 8} (2004), 180--211], the author gave a combinatorial algorithm for computing the Lusztig--Vogan bijection for $GL(n,\C)$.  However, that paper failed to mention one easy case that may sometimes arise, making the description of the algorithm incomplete.  This note fills in that gap.
\end{abstract}

\maketitle

\section{Introduction}
\label{sect:intro}

The main result of~\cite{a:phd,a:ekt} is the existence of a natural bijection between the set $\wts$ of dominant weights for $GL(n,\C)$ and the set $\ngo$ of pairs $(C,\cE)$, where $C$ is a nilpotent orbit, and $\cE$ is an irreducible equivariant vector bundle on $C$ (up to isomorphism).  Combinatorial objects called \emph{weight diagrams} play a key role, serving as an intermediary between $\wts$ and $\ngo$.  Specifically, there are maps
\[
\prj: \tD_n \to \ngo
\qquad\text{and}\qquad
\tpi: \tD_n \to \wts
\]
that both become bijections when restricted to the set $\sD_n \subset \tD_n$ of \emph{distinguished weight diagrams}.  The proofs of both assertions consist of giving algorithms for computing the inverse map.

The algorithm for the map $\ngo \to \sD_n$ involves choosing some (not necessarily distinguished) weight diagram $X \in \prj\inv(C,\cE)$, and then applying a sequences of ``moves'' to make $X$ distinguished.  At each step, one chooses a move to make by examining various properties, denoted by $\PP1(r)$, $\PP2(r)$, $\PP3(r)$, and $\PP4(r)$.

Unfortunately, the proofs in~\cite{a:ekt} do not cover the case where $\PP4(1)$ is false for the initial weight diagram.  Moreover, this omission is somewhat ``hidden,'' because the statements of various propositions purport to cover this case.  Happily, these mistakes are quite easy correct.  In this note, we indicate the appropriate corrections to statements in~\cite{a:ekt}, and we prove one new assertion, Proposition~7.9 below, to cover the missing case.

These mistakes came to light in November 2015, when David Vogan informed me of unexpected behaviour\footnote{This is a euphemism for ``The program crashed.''} in the software implementation~\cite{a:lvbgln} in a specific example.  I initially thought that this would be a matter of tracking down a software bug.  It turned out instead  to be a reasoning bug.  Indeed, the software faithfully implemented the algorithm of~\cite{a:ekt}, including the failure to handle the case where $\PP4(1)$ is false.  The software has now been updated to incorporate the corrections described below.  Vogan's example is reproduced in Section~\ref{sect:vogan}.

\section{Corrections to~\cite{a:ekt}}
\label{sect:corr}

\subsection{} On page~190, line~5, change
\begin{quote}
Every $Y$ is said to have properties $\pp1(1)$ and $\pp3(1)$, for
convenience.
\end{quote}
to
\begin{quote}
Every $Y$ is said to have properties $\pp1(1)$, $\pp2(1)$, and $\pp3(1)$, for
convenience.
\end{quote}
\subsection{} On page~190, lines~8--9, change
\begin{quote}
We suppose that the properties $\PP1(0)$, $\PP3(0)$, and
$\PP4(0)$ always hold, as do $\PP1(1)$ and $\PP3(1)$.
\end{quote}
to
\begin{quote}
We suppose that the properties $\PP1(0)$, $\PP2(0)$, $\PP3(0)$, and
$\PP4(0)$ always hold, as do $\PP1(1)$, $\PP2(1)$, and $\PP3(1)$.
\end{quote}
\subsection{} On page~197, in Table~2, the definitions of $q_5$ and $q_6$ define the parameter $r$ incorrectly.  They should both be amended to say
\[
r = -q_4(X) + 1.
\]
Also, the definition of $\tq_{5;ir}$ should be amended to consider the case $r = 1$ separately.  The corrected definition is
\begin{multline*}
\tq_{5;ir}(X) = \\
\begin{cases}
0 & \text{if $r = 1$,} \\
\max\{\prho X_{ir}-(\prho X_{i,r-1}+1),
  \prho X_{i,r-1}-\prho X_{ir}\}&\text{if $r$ is odd and $r \ne 1$,}\\
\max\{\prho X_{ir}-\prho X_{i,r-1},
  (\prho X_{i,r-1}-1)-\prho X_{ir}\}&\text{if $r$ is even.}\\
\end{cases}
\end{multline*}
\subsection{} On page~200, the statement of Proposition~7.8 should include ``$r > 1$'' as a hypothesis.  The corrected statement is:
\begin{prop78}
Suppose $r > 1$, and suppose rows $i$ and $i'$ of $X$ are such that \Br{} or \Cr{} might be
performed on them: in particular, they agree in their first $r-1$
entries, and intervening rows have length less than $r-1$.  Suppose
furthermore that $q_4(X) = -(r-1)$, and that row $i$ has length at
least $r$.
\begin{enumgen}{1}
\item Suppose $X$ has no entry at position $i'r$.  If $\tq_{5;ir} > 0$
and $\prho X_{ir} < \prho X_{i,r-1}$, $\Br X$ is well-behaved of order
$\ge 4$.\label{case:bcb}
\item Suppose that $X$ does have an entry at position $i'r$.  If
$X_{ir} < X_{i'r}$, then $\Cr X$ is well-behaved of order $\ge
4$.\label{case:bcc}
\end{enumgen}
\end{prop78}

\subsection{} On page~202, before the last sentence of \S7, insert the following new proposition:

\begin{prop79}
Suppose $q_4(X) = 0$, and let $r = 1$.  Let $i$ and $i' = i+1$ be two consecutive rows of $X$ such that $X_{ir} < X_{i'r}$.  Then $\Cr X$ is well-behaved of order~$\ge 4$.
\end{prop79}
\begin{proof}
As in the proof of Proposition~7.8, this move preserves $q_1$, $q_2$, and $q_3$.  The assumption implies that $X$ does not have $\PP4(1)$, and that $q_6(X) > 0$.  It is easy to see from the definition that
\[
\tq_{6;jr}(\Cr X) =
\begin{cases}
\tq_{6;jr}(X) & \text{if $j \ne i, i'$,} \\
\tq_{6,jr} - (X_{i'r} - X_{ir}) & \text{if $j = i$ or $i'$.}
\end{cases}
\]
Therefore, $q_6(\Cr X) < q_6(X)$.  If $q_6(\Cr X) = 0$, then $\Cr X$ has $\PP4(1)$, so $q_4(\Cr X) \le -1 < q_4(X)$, and hence $\Cr X$ is well-behaved of order~$4$.  On the other hand, if $q_6(\Cr X) > 0$, then we still have $q_4(\Cr X) = 0$.  In that case, we also have$q_5(\Cr X) = q_5(X) = 0$, so $\Cr X$ is well-behaved of order~$6$.
\end{proof}

\subsection{} On page~202, the last sentence of~\S7 should be amended to mention the new Proposition~7.9.  It should say:
\begin{quote}
The facts in Propositions~7.5, 7.7, 7.8 and 7.9 are collected and summarized in
Table~3.
\end{quote}
In addition, Table~3 should be amended to include the case covered by the new Proposition~7.9.  The corrected table is shown above.

\setcounter{table}{2}
\begin{table}
\begin{center}
\newcommand{\hem}{\hspace{1em}}
\begin{tabular}{|p{.5\textwidth}|p{.1\textwidth}|p{.28\textwidth}|}
\hline
{\em Conditions} & {\em Move} & {\em Well-Behavedness} \\
\hline
\hline

$q_4(X) \le -r$. & & Order $\ge 1$. \\
\hem $X_{i_ms}$ is lowerable, $X_{i_1r}$ is
raisable, and $\prho X_{i_ks} - \prho X_{i_kr} \ge 1$ for $k = 1,
\ldots, m$. &
\Ar & \\
\hem $X_{i_1s}$ is raisable, $X_{i_mr}$ is
lowerable, and $\prho X_{i_ks} - \prho X_{i_kr} \le -1$ for $k = 1,
\ldots, m$. &
\Al & \\

\hline
$m = 1$; $q_4(X) = -(r - 1)$.  $X_{i_1s}$ is lowerable. $\tilde
q_{3;i_1r}(X) \ne 0$. If $X_{jr} = X_{ir}$, then $j \ge i$ (for \Ar)
or $j \le i$ (for \Al).& & Order $= 1$. \\
\hem $\prho X_{i_1s} - \prho X_{i_1r} \ge 2$.\par
\hem $\prho X_{i_1s} - \prho X_{i_1r} \le -2$. &
\Ar \par \Al & \\
\hem \hrulefill &\hrulefill &\\
\hem $r$ odd, $s$ even, and $\prho X_{i_1s} - \prho X_{i_1r} = 1$.\par
\hem $r$ even, $s$ odd, and $\prho X_{i_1s} - \prho X_{i_1r} = -1$. &
\Ar \par \Al & \\

\hline
$q_4(X) = -(r - 1)$; $\tq_{5;i_1r}(X) > 0$\par
\hem $\prho X_{i_1,r-1} > \prho X_{i_1r}$.\par
\hem $\prho X_{i_1,r-1} < \prho X_{i_1r}$.\par
\hem \hrulefill \par
\hem $X_{ir} < X_{i'r}$. &
\ \par \Br \par \Bl \par \hrulefill \par \Cr &
Order $\ge 4$. \\

\hline
$q_4(X) = 0$; $i' = i+1$; $r = 1$\par
\hem $X_{ir} < X_{i'r}$. &
\ \par
\Cr &
Order $\ge 4$. \\

\hline

\end{tabular}
\end{center}
\caption{Well-behavedness of moves under various
hypotheses\label{tbl:behavior}}
\end{table}

\subsection{} On page~203, the statements of Proposition~8.3 and Corollary~8.4 should include ``$r > 1$'' as a hypothesis. The corrected statements are:
\begin{prop83}
Suppose $r > 1$. If $\prho X$ has $\pp3(r)$ and $\pp4(r-1)$, then it also has
$\pp4(r)$.
\end{prop83}

\begin{cor84}
Suppose $r > 1$, and that $\prho X$ has $\PP1(r-1)$, $\PP2(r-1)$, $\PP3(r-1)$, and
$\PP4(r-1)$.  If it does not have $\PP4(r)$, then it also does not
have $\PP3(r)$. \qed
\end{cor84}

\subsection{} On page~204, the statement of Proposition~8.5 should be amended to say

\begin{prop85}
Suppose $\prho X$ has $\PP1(r-1)$, $\PP2(r-1)$, $\PP3(r-1)$, and
$\PP4(r-1)$.  If it does not have $\PP3(r)$ or $\PP4(r)$, then $X$ satisfies some
hypothesis in the left-hand column of Table \ref{tbl:behavior}.
\end{prop85}

Its proof should have the following new paragraph inserted at the beginning:

\begin{proof}
Suppose first that $r = 1$.  In this case, $\PP3(1)$ holds automatically, so it must fail to have $\pp4(1)$.  This means that there is some entry $X_{i1}$ that is smaller than its column-successor $X_{i'1}$, where necessarily $i' = i+1$.  In other words, we are in the setting of Proposition~7.9, and we can do move \Cr.  

For the remainder of the proof, assume that $r > 1$. Corollary~8.4 tells us that $\PP3(r)$ must fail, so we know that $q_5(X) > 0$.  [The rest of the proof is as in~\cite{a:ekt}.]
\end{proof}

\section{Vogan's Example}
\label{sect:vogan}

We will compute the distinguished weight diagram in $\prj\inv(C,\cE)$ for the pair $(C,\cE) = ([3^3,2^2];((5,0,-3),(4,-6)))$.
In the calculation below, which follows the pattern of~\cite[\S9.3]{a:ekt}, Steps~2 and~3  rely on the new Proposition~7.9.  Note that for the initial weight diagram $X$, $\prho X$ fails to satisfy $\PP4(1)$.  Moreover, it only satisfies the hypothesis of the newly added row in Table~3.  Thus, this example demonstrates the incompleteness of the description of the algorithm given in~\cite{a:ekt}.

\begin{enumerate}
\item Choose a weight diagram $X$ such that $\prj(X) = ([3^3,2^2];((5,0,-3),(4,-6)))$.
\[
X = \begin{tableau}
:.5.0.0\\:.0.0.0\\:.{-3}.0.0\\:.4.0\\:.{-6}.0\\
\end{tableau}
\qquad\qquad\qquad\qquad
\prho X = \begin{tableau}
:.9.4.2\\:.0.2.0\\:.{-5}.0.{-2}\\:.6.{-2}\\:.{-10}.{-4}\\
\end{tableau}
\]

\item Perform \Cr{} with $r = 1$, $i = 3$, and $i' = 4$.
\[
X = \begin{tableau}
:.5.0.0\\:.0.0.0\\:.4.0\\:.{-3}.0.0\\:.{-6}.0\\
\end{tableau}
\qquad\qquad\qquad\qquad
\prho X = \begin{tableau}
:.9.4.2\\:.0.2.0\\:.6.0\\:.{-5}.{-2}.{-2}\\:.{-10}.{-4}\\
\end{tableau}
\]

\item Perform \Cr{} with $r = 1$, $i = 2$, and $i' = 3$.
\[
X = \begin{tableau}
:.5.0.0\\:.4.0\\:.0.0.0\\:.{-3}.0.0\\:.{-6}.0\\
\end{tableau}
\qquad\qquad\qquad\qquad
\prho X = \begin{tableau}
:.9.4.2\\:.6.2\\:.0.0.0\\:.{-5}.{-2}.{-2}\\:.{-10}.{-4}\\
\end{tableau}
\]

\item Perform \Ar{} four times with $s = 1$, $r = 2$, on rows $1$ and $2$.
\[
X = \begin{tableau}
:.3.2.0\\:.2.2\\:.0.0.0\\:.{-3}.0.0\\:.{-6}.0\\
\end{tableau}
\qquad\qquad\qquad\qquad
\prho X = \begin{tableau}
:.7.6.2\\:.4.4\\:.0.0.0\\:.{-5}.{-2}.{-2}\\:.{-10}.{-4}\\
\end{tableau}
\]

\item Perform \Al{} five times with $s = 1$, $r = 2$, on rows $4$ and $5$.
\[
X = \begin{tableau}
:.3.2.0\\:.2.2\\:.0.0.0\\:.{-1}.{-2}.0\\:.{-3}.{-3}\\
\end{tableau}
\qquad\qquad\qquad\qquad
\prho X = \begin{tableau}
:.7.6.2\\:.4.4\\:.0.0.0\\:.{-3}.{-4}.{-2}\\:.{-7}.{-7}\\
\end{tableau}
\]

\item Perform \Ar{} with $s = 1$, $r = 3$ on row $1$.
\[
X = \begin{tableau}
:.2.2.1\\:.2.2\\:.0.0.0\\:.{-1}.{-2}.0\\:.{-3}.{-3}\\
\end{tableau}
\qquad\qquad\qquad\qquad
\prho X = \begin{tableau}
:.6.6.3\\:.4.4\\:.0.0.0\\:.{-3}.{-4}.{-2}\\:.{-7}.{-7}\\
\end{tableau}
\]

\item Perform \Al{} with $s=2$, $r = 3$ on row $4$.
\[
X = \begin{tableau}
:.2.2.1\\:.2.2\\:.0.0.0\\:.{-1}.{-1}.{-1}\\:.{-3}.{-3}\\
\end{tableau}
\qquad\qquad\qquad\qquad
\prho X = \begin{tableau}
:.6.6.3\\:.4.4\\:.0.0.0\\:.{-3}.{-3}.{-3}\\:.{-7}.{-7}\\
\end{tableau}
\]

\item Perform \Br{} on rows $1$ and $2$.
\[
X = \begin{tableau}
:.2.2\\:.2.2.1\\:.0.0.0\\:.{-1}.{-1}.{-1}\\:.{-3}.{-3}\\
\end{tableau}
\qquad\qquad\qquad\qquad
\prho X = \begin{tableau}
:.6.6\\:.4.4.3\\:.0.0.0\\:.{-3}.{-3}.{-3}\\:.{-7}.{-7}\\
\end{tableau}
\]

\item Perform \Ar{} with $s = 2$, $r = 3$ on row $2$.
\[
X = \begin{tableau}
:.2.2\\:.2.1.2\\:.0.0.0\\:.{-1}.{-1}.{-1}\\:.{-3}.{-3}\\
\end{tableau}
\qquad\qquad\qquad\qquad
\prho X = \begin{tableau}
:.6.6\\:.4.3.4\\:.0.0.0\\:.{-3}.{-3}.{-3}\\:.{-7}.{-7}\\
\end{tableau}
\]

\end{enumerate}

From the last step, we obtain
\[
\gamma([3^3,2^2];((5,0,-3),(4,-6))) = \tpi(X) =
(6,6,4,4,3,0,0,0,-3,-3,-7,-7).
\]


\begin{thebibliography}{0}

\bibitem{a:phd}
P.~Achar, {\em Equivariant coherent sheaves on the nilpotent cone for complex
  reductive Lie groups}, Ph.D. thesis, Massachusetts Institute of Technology,
  2001.

\bibitem{a:ekt}
P.~Achar, {\em On the equivariant $K$-theory of the nilpotent cone in the
  general linear group}, Represent. Theory {\bf 8} (2004), 180--211.

\bibitem{a:lvbgln}
P.~Achar, {\em The Lusztig--Vogan bijection for $GL_n$}, software available for
  download from {\tt http://www.math.lsu.edu/\char"7E pramod/}, 2004.

\end{thebibliography}
\end{document}